\documentclass[11pt]{article}
\usepackage[dvipdfmx]{graphicx,epsfig}  
\usepackage{subfigure}
\usepackage{epstopdf}
\usepackage{amssymb,amsmath}
\usepackage{caption}
\usepackage{float}
\usepackage{secdot}
\usepackage{amsthm}

\newtheorem{thm}{Theorem}[section]
\newtheorem{lem}[thm]{Lemma}

\newtheorem{prop}[thm]{Proposition}
\newtheorem{cor}[thm]{Corollary}

\newtheorem{defn}[thm]{Definition}
\newtheorem{rmk}[thm]{Remark}
\newtheorem{exmp}[thm]{Example}

\usepackage{amssymb}

\begin{document}
\makeatletter{\renewcommand*{\@makefnmark}{}
\footnotetext{2010 Mathematics Subject Classification. Primary  57M25 ; Secondary 57M27.}\makeatother}

\begin{center}
\Large{AN UNKNOTTING INDEX FOR VIRTUAL KNOTS} 
\end{center} 
\vspace{.5cm}
\begin{center}
K. KAUR, 
S. KAMADA\footnote{Corresponding author.  Supported by JSPS KAKENHI Grant Numbers 26287013 and 15F15319.}, 
A. KAWAUCHI\footnote{{Supported by JSPS KAKENHI Grant Number 24244005.}} and 
M. PRABHAKAR

\end{center}

\begin{center} \textbf{Abstract} \end{center}

In this paper we introduce the notion of an unknotting index for virtual knots. 
We give some examples of computation by using writhe invariants, and 
discuss a relationship between the unknotting index and the virtual knot module. In particular, we show that for any non-negative integer $n$ there exists a virtual knot whose unknotting index is $(1,n)$. 


\section{Introduction}
\label{}

\begin{figure}[!ht] 
\centering  
\includegraphics[scale=.4]{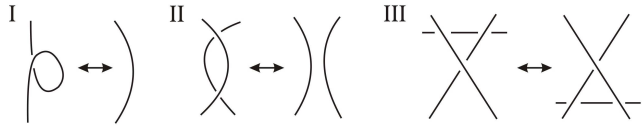} 
\caption{Classical Reidemeister moves: RI, RII and RIII}
\label{fig:classical_reimet}
\end{figure}

\begin{figure}[!ht] 
\centering  
\includegraphics[scale=.4]{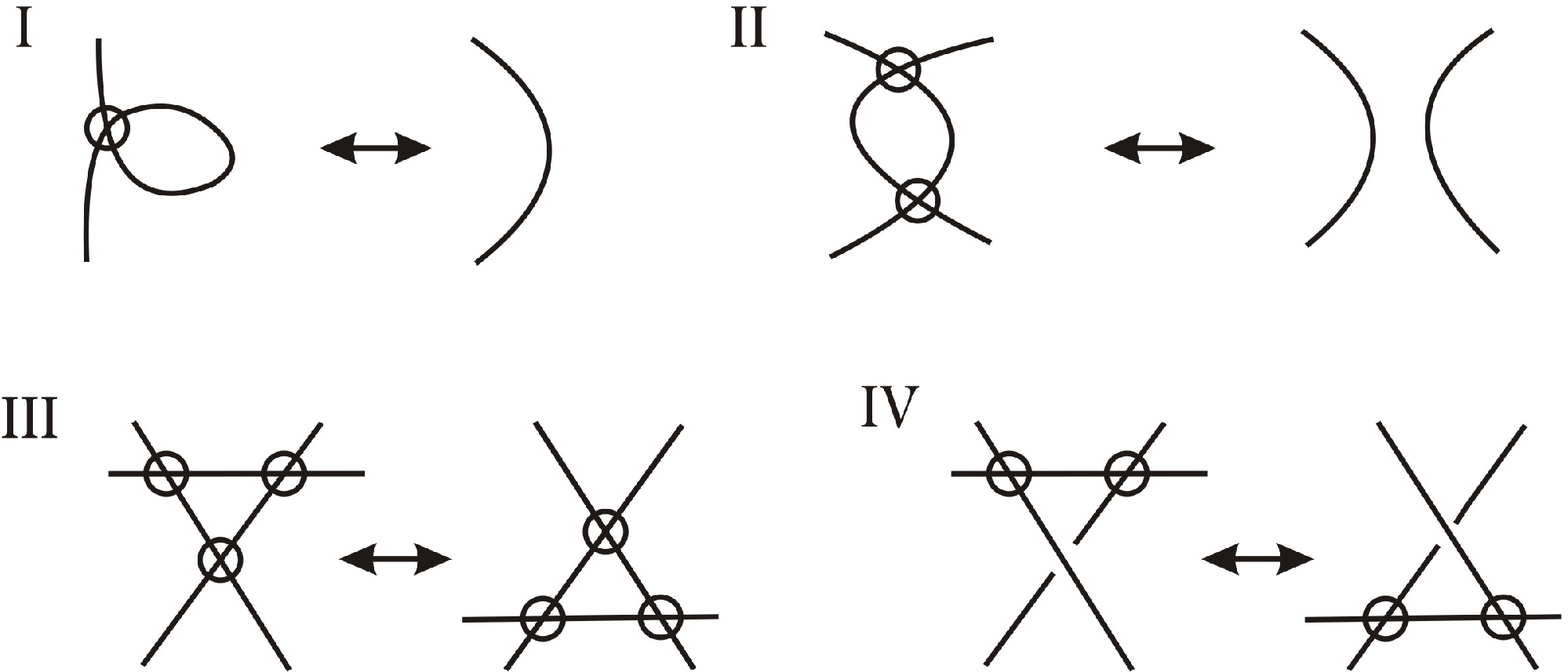} 
\caption{Virtual Reidemeister moves: VRI, VRII, VRIII and VRIV}
\label{fig:virtual_reimet}
\end{figure}

Virtual knot theory was introduced by L.~H.~Kauffman \cite{kauffman1998virtual} as a generalization of classical knot theory.  
A diagram of a virtual knot may have virtual crossings which are encircled by a small circle and are not regarded as (real) crossings.  
A virtual knot is an equivalence class of diagrams where the equivalence is generated by moves in Figs.~\ref{fig:classical_reimet} and \ref{fig:virtual_reimet}. 

 The purpose of this paper is to introduce an unknotting index $U(K)$ for a virtual knot, 
  whose idea is an extension of the usual unknotting  number for classical knots.  
  The unknotting index of a virtual knot is a pair of non-negative integers which is considered as an ordinal number with respect to the dictionary order.  
The definition is given in Section~\ref{sect:prelimnaries}.  
As the usual unknotting number for classical knots, it is not easy to determine 
the value of the unknotting index of a given virtual knot.  
We provide  some examples of computation by using writhe invariants of virtual knots in Section~\ref{sect:main_writhe}. 
In Section~\ref{sect:main_module} 
we discuss the unknotting index using the virtual knot module. An estimation of the unknotting index using the index $e(M)$ of the virtual knot module $M$ is given (Theorem~\ref{thm:module}).  Then we show that 
for any non-negative integer $n$, there exists a virtual knot $K$ with $U(K) = (0,n)$ 
(Example~\ref{example:U0n} or Proposition~\ref{thm:U0n}) 
and there exists a virtual knot $K$ with $U(K) = (1,n)$ (Theorem~\ref{thm:U1n}).

\section{An unknotting index for virtual knots}
\label{sect:prelimnaries} 

In this section, we define an unknotting index for virtual knots.  

Let $D$ be a diagram of a virtual knot $K$.  
For a pair $(m,n)$ of non-negative integers,  
we say that $D$ is {\it $(m,n)$-unknottable} if the diagram $D$ is changed into a diagram of the trivial knot by changing  $m$ crossings of $D$ into virtual crossings and by applying crossing change operations on $n$ crossings of $D$.   
Let ${\mathcal U}(D)$ denote the set of such pairs $(m,n)$.  
Note that $D$ is $(c(D), 0)$-unknottable, where $c(D)$ is the number of crossings of $D$. 
Thus, ${\mathcal U}(D)$ is non-empty.  

\begin{defn}{\rm 
\begin{itemize}
\item 
The {\it unknotting index}  of a diagram $D$, denoted by $U(D)$, is the minimum among all pairs $(m,n)$ such that $D$ is $(m,n)$-unknottable.  Namely, $U(D)$ is the minimum among 
 the family ${\mathcal U}(D)$.    (The minimality is taken with respect to the dictionary order.) 
\item 
The {\it unknotting index} of a virtual knot $ K$, denoted by $U(K)$, is the minimum  
among all pairs $(m,n)$ such that $K$ has a diagram $D$ which is $(m,n)$-unknottable. Namely, 
$U(K)$ is the minimum among the family 
$\{ U(D) \mid \mbox{$D$ presents $K$} \}$.  
\end{itemize}
}\end{defn}

\begin{figure}[!ht] 
\centering
\includegraphics[scale=0.35]{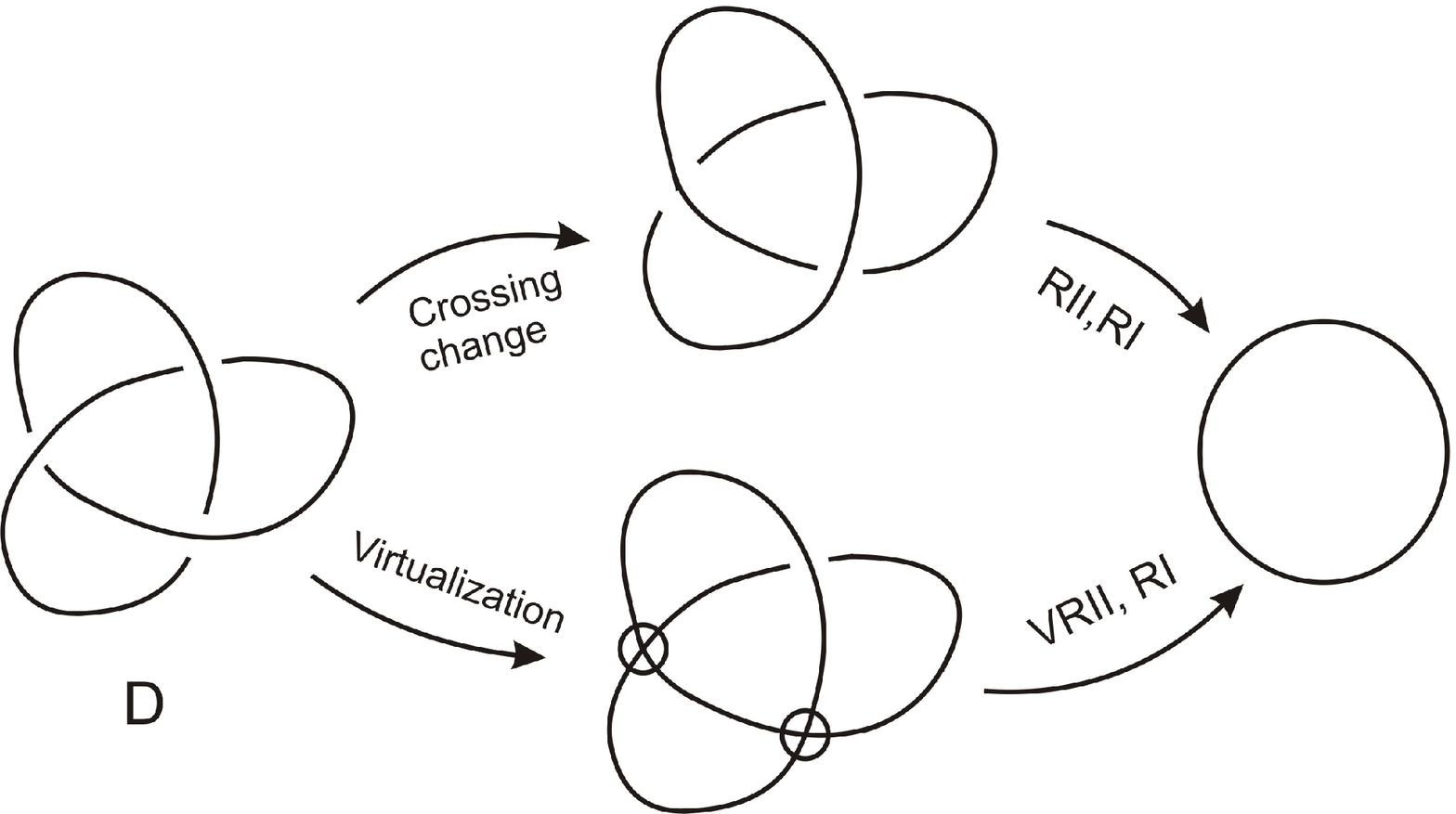}  
 \caption{}
 \label{example_of_generalized_unknotting_number}

\end{figure}

For example, when $D$ is the diagram on the left of  
Fig.~\ref{example_of_generalized_unknotting_number}, the figure shows that $(0,1)$ and $(2,0)$ belong to ${\mathcal U}(D)$.  Since $(0,0)$ is not an element of ${\mathcal U}(D)$, we have $U(D)= (0,1)$.  (Actually, ${\mathcal U}(D) =\{ (0,1), (0,2), (1,1), (2,0), (2,1), (3,0)\}$. In order to determine $U(D)$ 
 it is often unnecessary to find all elements of ${\mathcal U}(D)$.) 
Thus, for the left-handed trefoil $K$, $U(K) = (0,1)$. 
(Note that a virtual knot $K$ is trivial iff $U(K)=(0,0)$.)   

\begin{figure}[!ht] 
\centering
\includegraphics[scale=0.35]{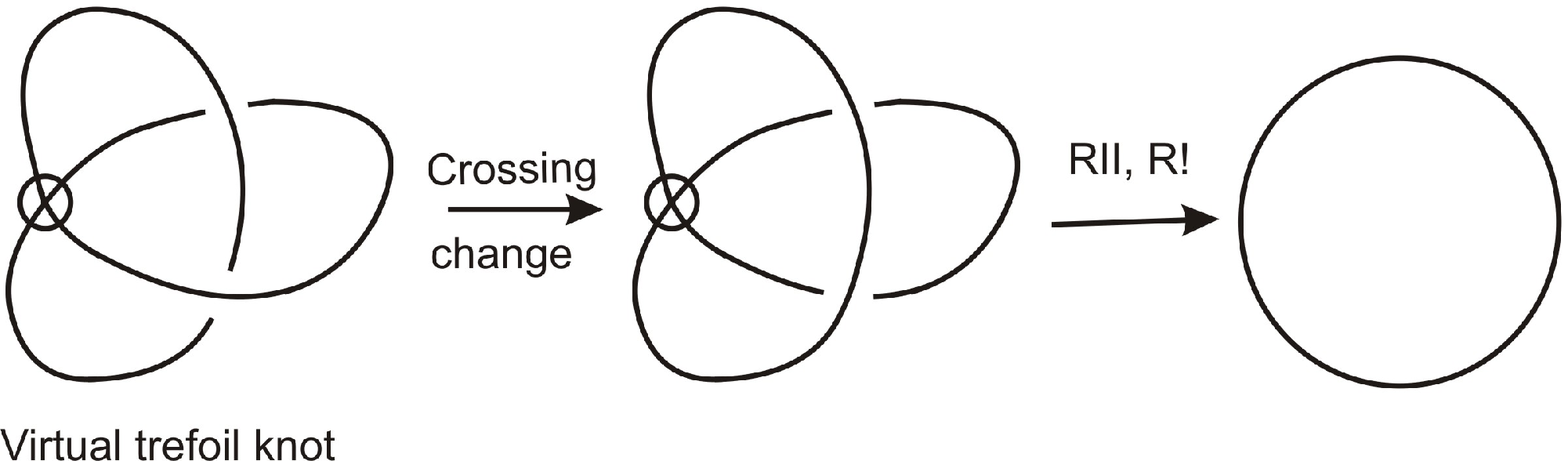} 
 \caption{}
 \label{fig:virtual_trefoil}
\end{figure}

For the (left-handed) virtual trefoil, which is presented by the diagram of the left of  
Fig.~\ref{fig:virtual_trefoil},   $U(K)=(0,1)$.

A {\it flat virtual knot diagram} is a virtual knot diagram by forgetting the over/under-information of every real crossing.  
A {\it flat virtual knot} is an equivalence class of flat virtual knot diagrams by 
{\it flat Reidemeister moves} which are Reidemeister moves (Figs.~\ref{fig:classical_reimet} and \ref{fig:virtual_reimet}) without the over/under-information.   

\begin{lem}\label{lem:flat}
Let $K$ be a virtual knot and $\overline{K}$ its flat projection. If $\overline{K}$ is non-trivial, then 
$U(K) \geq (1,0)$. 
\end{lem}

\noindent \textbf{Proof:} 
Assume that $U(K) = (0,n)$ for some $n$. There is a diagram $D$ of $K$ such that $D$ becomes a diagram of the trivial knot by $n$ crossing changes.  The flat virtual knot diagram $\overline{D}$ obtained from $D$ is a diagram of the trivial knot, which implies $\overline{K}$ is trivial.  \hfill $ \square $

\begin{figure}[!ht] 
\centering
\includegraphics[scale=0.4]{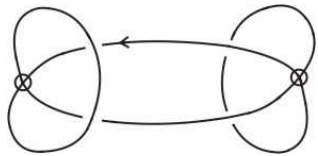}
\caption{Kishino's knot}
\label{kishhh}
\end{figure}

\begin{exmp}{\rm 
Let $D$ be the diagram depicted in Fig.~\ref{kishhh} and $K$ the virtual knot presented by $D$, which is  
called {\it Kishino's knot}.  
It is nontrivial as a virtual knot 
(cf. \cite{BartholomewFenn, FennTuraev, Kadokami, KishinoSatoh, miyazawa2008}) and the flat projection $\overline{K}$ is nontrivial as a flat virtual knot (cf. \cite{FennTuraev, Kadokami}).  
Thus $U(K) \geq (1,0)$.  It is easily seen that 
$D$ is $(1,0)$-unknottable.  Therefore, $U(K) =(1,0)$.  
}\end{exmp}

\section{Gauss diagrams and writhe invariants}
\label{sect:some_inv}

We recall Gauss diagrams and writhe invariants of virtual knots.  
In what follows, we assume that virtual knots are oriented.

A Gauss diagram of a virtual knot diagram $D$ is an oriented circle where the pre-images of the over crossing and under crossing of each crossing are connected by a chord. To indicate the over/under-information,  chords are directed from the over crossing to the under. Corresponding to each crossing $c$ in a virtual knot diagram, there are two points $ \overline{c}$ and $ \underline{c} $ which present the over crossing  and the under crossing  for $c$ in the oriented circle. The sign of each chord is the sign of the corresponding crossing.  The sign of a crossing or a chord $c$ is also called the {\it writhe} of $c$ and denoted by $w(c)$ in this paper.  

\begin{figure}[!ht] 
\begin{center}  
\includegraphics[scale=0.8]{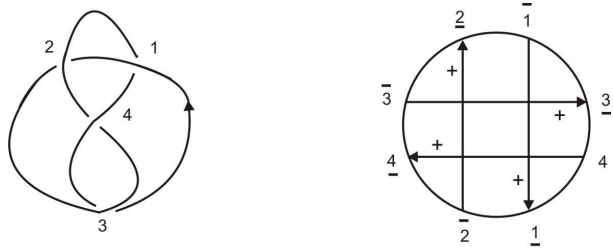} 
\caption{Gauss diagram for a diagram of the figure eight knot}
\label{fig:Gauss}
\end{center} 
\end{figure}

In terms of a Gauss diagram, virtualizing a crossing of a diagram corresponds to elimination of a chord, and a crossing change corresponds to changing the direction and the sign of a chord. 
In what follows, for a Gauss diagram, we refer to changing the direction and the sign of a chord $c$ as a {\it crossing change} at $c$. 

\vspace{0.3cm}

Writhe invariants were defined by some researchers independently.  
The writhe polynomial was defined by Cheng--Gau \cite{cheng2013polynomial}, which  
is equivalent to the affine index polynomial defined by Kauffman~\cite{kauffman2013}.  
Satoh--Taniguchi~\cite{SatohTaniguchi2014} introduced the $k$-th writhe $J_k$ for each $k \in {\mathbb Z}\setminus \{0\}$, which is indeed a coefficient of the affine index polynomial.    
(Invariants related to these  
are found in Cheng~\cite{cheng2014polynomial}, 
Dye~\cite{Dye2013} and Im--Kim--Lee~\cite{ImKimLee2013}.)   

The writhe polynomial $W_G(t)$ stated below is the one in \cite{cheng2013polynomial} multiplied by $t^{-1}$.  This convention makes it easier to see the relationship between  $W_G(t)$ and     
the affine index polynomial or the $k$-th writhe $J_k$. 

\vspace{0.3cm}

Let $c$ be a chord of a Gauss diagram $G$.  
Let $r_{+}$  (respectively, $r_{-}$) be the number of positive  (resp. negative) chords intersecting with $c$ transversely from right to left when we see them from the tail toward the head of $c$.   
Let $l_{+}$  (respectively, $l_{-}$) be the number of positive (resp. negative) chords intersecting with $c$ transversely from left to right.   
The {\it  index}  of $c$ is defined as 
$${\rm Ind}(c)=r_{+}-r_{-}-l_{+}+ l_{-}.$$ 
The {\it writhe polynomial} $W_{G}(t)$ of $G$ is defined by 
$$W_{G}(t) 
=   \sum_{c ~:~ {\rm Ind}(c) \neq 0}     w(c) t^{{\rm Ind}(c)}. $$ 

For each integer $k$, the {\it $k$-th writhe} $J_k(G)$ is the number of positive chords with index $k$ minus that of negative ones with index $k$.  Then 
$$W_{G}(t) 
=   \sum_{k \neq 0}     J_k(G) t^{k}. $$ 
The  $k$-th writhe $J_k(G)$ is an invariant of the virtual knot presented by $G$ when $k \neq 0$.  

The {\it writhe polynomial} $W_{K}(t)$ and the 
{\it $k$-th writhe} $J_k(K)$
of a virtual knot $K$  is defined by 
those for a Gauss diagram presenting $K$.  

The odd writhe $J(G)$ defined by Kauffman \cite{kauffman2004} is  
$$ \sum_{c \in Odd(G)} w(c),$$ 
where $w(c)$ denotes the writhe and $Odd(G)$ is the set of chords with odd indices.

\section{Unknotting indices of some virtual knots}
\label{sect:main_writhe}

We give examples of computation of the unknotting indices for some virtual knots using writhe invariants.  

\begin{lem}[\cite{SatohTaniguchi2014}]
\label{lem:cc} 
Let $G$ and $G'$ be Gauss diagrams such that $G'$ is obtained from $G$ by a crossing change.  
Then one of the following occurs. 
\begin{itemize}
\item[$(1)$] 
$ W_G(t) - W_{G'}(t) = \epsilon ( t^m + t^{-m})$ for some integer $m \in {\mathbb Z}\setminus \{0\}$ and $\epsilon \in \{ \pm 1\}$.  
\item[$(2)$] 
$ W_G(t) = W_{G'}(t)$ 
\end{itemize} 
\end{lem}

\noindent \textbf{Proof:} 
Let $c$ be the chord of $G$ such that a crossing change at $c$ changes $G$ to $G'$. Let $c'$ be the chord of $G'$ obtained from $c$ by the crossing change.   
Then ${\rm Ind}(c') = -{\rm Ind}(c)$ and $w(c') = -w(c)$.  For any chord $d$ of $G$, except $c$,  the index ${\rm Ind}(d)$ and the writhe $w(d)$ are preserved by the crossing change at $c$.  When ${\rm Ind}(c) \neq 0$, let $m= {\rm Ind}(c)$ and $\epsilon = w(c)$, and we obtain the first case. 
When ${\rm Ind}(c) = 0$, we have the second case. \hfill $ \square $ 

\vspace{0.3cm} 

By this lemma, we have the following. 

\begin{prop}[cf. Theorem~1.5 of \cite{SatohTaniguchi2014}]
\label{prop:boundsUbyS}
Let $K$ be a virtual knot. 
\begin{itemize} 
\item[$(1)$]  If $J_k(K) \neq J_{-k}(K)$ for some $k \in {\mathbb Z}\setminus \{0\}$, then $(1,0) \leq U(K)$. 
\item[$(2)$]  $(0,   \sum_{k \neq 0} | J_k(K) |  /2 ) \leq U(K)$.  
\end{itemize}
\end{prop}

\noindent \textbf{Proof:} 
(1) Suppose that $(1,0) > U(K)$. This implies that 
$K$ has a Gauss diagram $G$ which can be transformed  by crossing changes into a Gauss diagram presenting the unknot. 
By Lemma~\ref{lem:cc}, the writhe polynomial $W_G(t)$ must be reciprocal, i.e., $J_k(G) = J_{-k}(G)$ for all $k \in {\mathbb Z}\setminus \{0\}$. 
(This is (i) of Theorem~1.5 of \cite{SatohTaniguchi2014}.) 
This contradicts the hypothesis.  Thus, $(1,0) \leq U(K)$. 

(2) If $(1,0) \leq U(K)$, then the inequality holds. Thus, we consider the case of $(1,0) > U(K) =(0, n)$.  
Let $G$ be a Gauss diagram of $K$ which can be transformed  by $n$ crossing changes into a Gauss diagram presenting the unknot.  
By Lemma~\ref{lem:cc}, we see that $ \sum_{k\neq 0} |J_k(K) | \leq 2n$.  \hfill $ \square $ 

\vspace{0.3cm} 

\begin{rmk}{\rm 
If $K$ has a Gauss diagram $G$ which can be transformed  by crossing changes into a Gauss diagram presenting the unknot, then 
$J_k(G) = J_{-k}(G)$ for all $k \in {\mathbb Z}\setminus \{0\}$.  In this case, 
$ \sum_{k\neq 0} |J_k(K) | /2 =  \sum_{k >  0} |J_k(K) | =  \sum_{k < 0} |J_k(K) |$.  (This is (ii) of Theorem~1.5 of \cite{SatohTaniguchi2014}.) 
Thus, we may say that $(0,   \sum_{k > 0} | J_k(K) | ) \leq U(K)$ instead of the inequality in the second assertion of Proposition~\ref{prop:boundsUbyS}.   
}\end{rmk}

\begin{cor}
\label{boundsUbyA}
Let $K$ be a virtual knot and $J(K)$ the odd writhe.  Then 
$$ (0,  | J(K)  / 2 | ) \leq  U(K).$$
\end{cor}

\noindent \textbf{Proof:} 
By definition, $| J(K)  | \leq \sum_{k \neq 0} | J_k(K) |$.  The inequality in the second assertion of Proposition~\ref{prop:boundsUbyS} implies the desired inequality.   \hfill $ \square $ 

\vspace{0.3cm} 

We  consider unknotting indices for some virtual knots.  

\begin{figure}[!ht] 
{
\centering
\subfigure[] 
{
\includegraphics[scale=0.34]{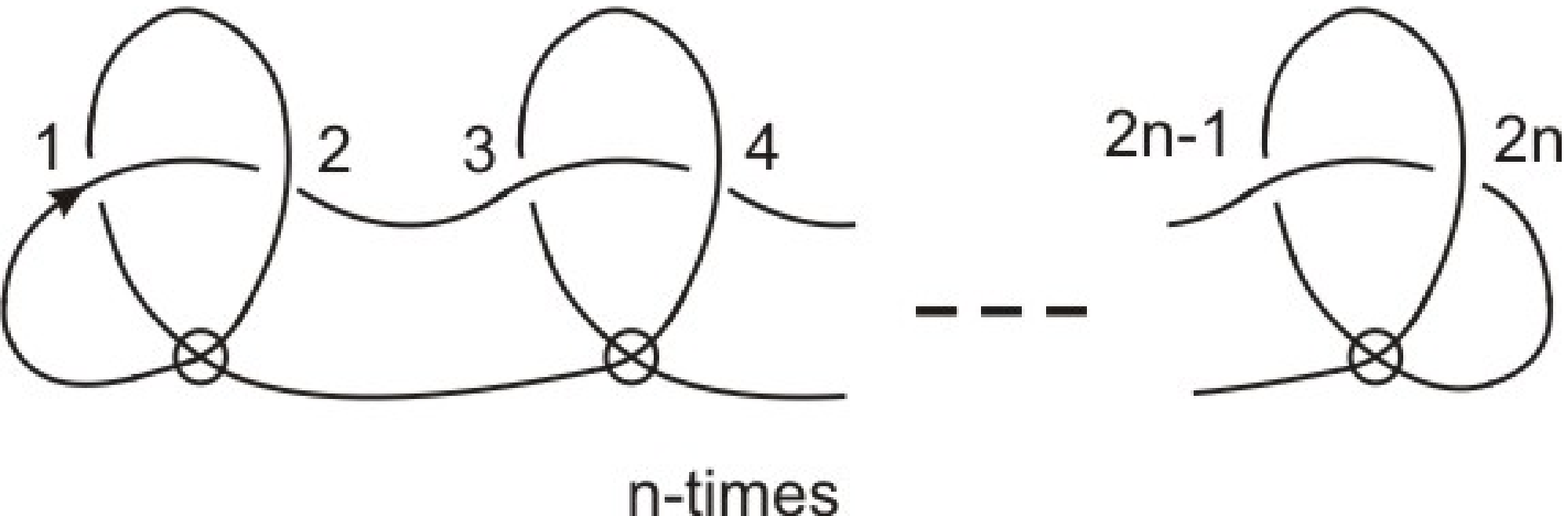}
\label{composition_trefoil_1}
}
\hspace{.4cm}
\subfigure[  ]
{
\includegraphics[scale=0.34]{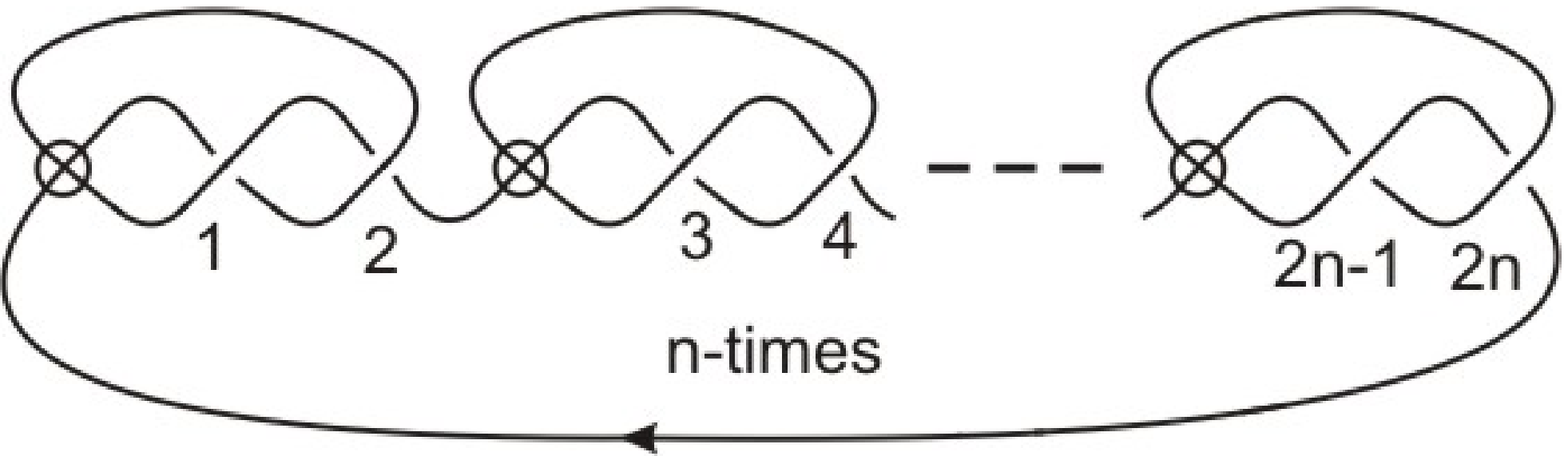} 
\label{composition_trefoil_2}
 }
\caption{} 
}
\end{figure}

\begin{exmp}\label{example:U0n}
{\rm Let $K$ be the virtual knot presented by a diagram $D$ depicted in Fig.~\ref{composition_trefoil_1} or Fig.~\ref{composition_trefoil_2}.  Then $U(K)=(0,n)$. 
} \end{exmp}

\noindent \textbf{Proof:} 
Since all crossings in $D $ are odd and have the same writhe, the odd writhe of $K$ is $ \pm 2n $. By Corollary~\ref{boundsUbyA}, $ (0, n) \leq U(K) $. By crossing changes at crossings labelled with even numbers, we obtain a diagram of the trivial knot. Therefore $ (0,n) \geq U(K) $. 
 \hfill $ \square $ 

\vspace{0.3cm} 

Let $p$ be an odd integer.  
By a {\it standard diagram} of the $(2,p)$-torus knot, we mean a diagram obtained from a diagram of the $2$-braid $\sigma_1^p$ by taking closure.  

\begin{exmp}
If $K$ is a virtual knot presented by a diagram $D$ which is obtained from a standard diagram of the $(2,p)$-torus knot by virtualizing some crossings, then  
$$ U(K) \leq ( 0, c/2  ) , $$ 
where $c$ is the number of crossings of $D$.  Moreover, if $c$ is even then $ U(K) = (0, c/2 )$.
\end{exmp}

\begin{figure}[!ht]
\begin{center}
\includegraphics[scale=0.38]{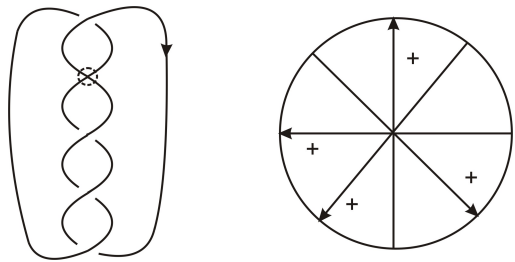}
\caption{ Closure of a virtual 2-braid with $U(K)=(0,2)$}
\label{example of complete virtual knot}
\end{center}
\end{figure}

\begin{figure}[!ht]
\begin{center}
\includegraphics[scale=.5]{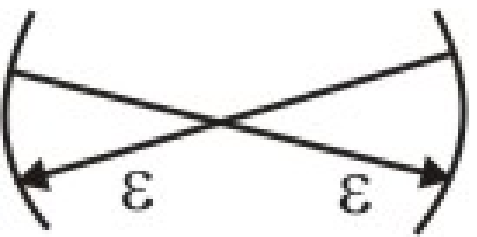}
\caption{ }
\label{fig:t1}
\end{center}
\end{figure}

\noindent \textbf{Proof:} 
Let $G$ be the Gauss diagram corresponding to the diagram $D$ of $K$.  An example is shown in Fig.~\ref{example of complete virtual knot}. When $ c\geq 2$, there exists a pair of chords as in Fig.~\ref{fig:t1}, which can be removed by a crossing change at one of the chords and 
an RII  move. If the resulting Gauss diagram still has such a pair of chords, then we remove the chords by a crossing change and an RII move. In this way, we can change $G$ into a Gauss diagram presenting the trivial knot by at most $c/2$ crossing changes. Hence, $ U(K)\leq (0, c/2)$. 
If $c$ is even, then each chord of $G$ is an odd chord.  The odd writhe of $K$ is $\pm c$.  
By Corollary~\ref{boundsUbyA},  $(0, c/2) \leq U(K)$. Hence $ U(K)= (0, c/2 )$.  \hfill $ \square $ 


Next, we consider a virtual knot presented by a diagram obtained from a diagram of a twisted knot.  

\begin{figure}[!ht] 
{
\centering
\subfigure[] 
{
\includegraphics[scale=0.35]{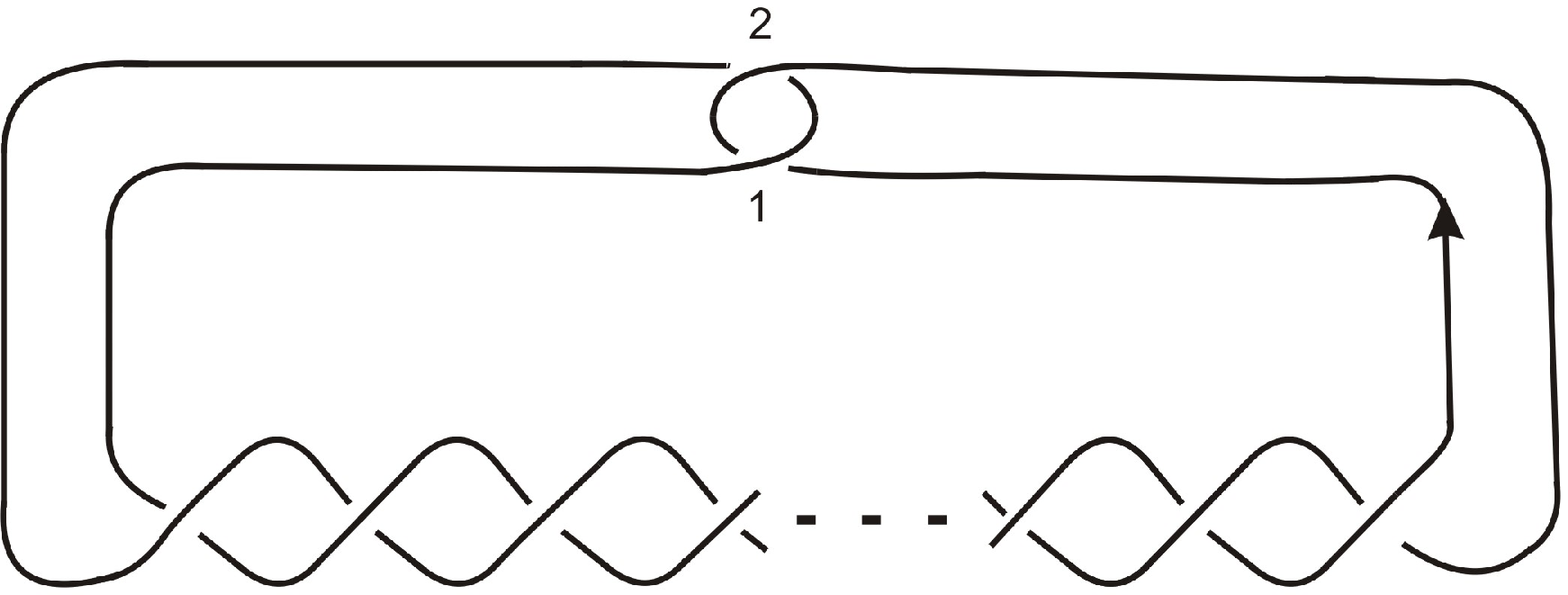}
\label{Twisted_Knot_Diagram}
}
\hspace{.1cm}
\subfigure[]
{
\includegraphics[scale=0.38]{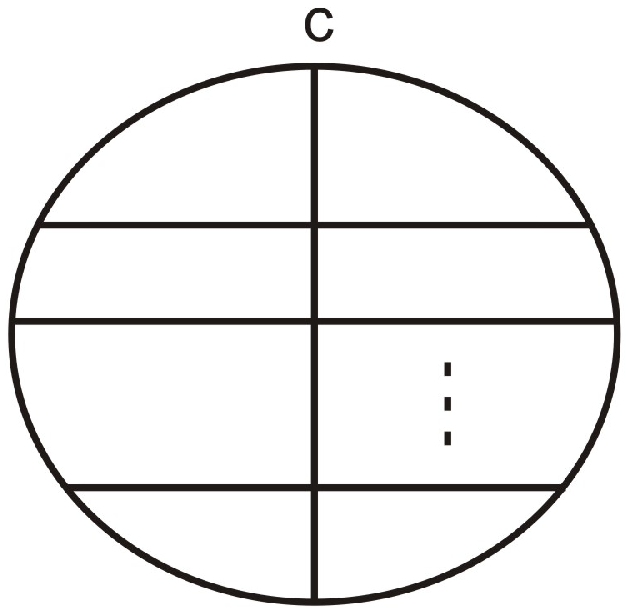} 
\label{c_chord}
}
 \caption{}
}
\end{figure}

Let $ D_0 $ be the diagram illustrated in Fig.~\ref{Twisted_Knot_Diagram}, which we call a 
{\it standard diagram} of a twisted knot.  All crossings of $D_0$ except the crossings labeled $1$ and $2$ are positive.  The signs of crossings labeled $1$ and $2$ are the same, which depends on the number of crossings.

Let $D$ be a diagram obtained from $D_0$ by virtualizing some crossings.  Let $K$ be the virtual knot presented by $D$. 

If the crossings labeled 1 and 2 in Fig.~\ref{Twisted_Knot_Diagram} are both virtualized in $D$, then the diagram $D$ presents the trivial knot. Hence $U(K)=(0,0)$.  

If the crossings labeled 1 and 2 in Fig.~\ref{Twisted_Knot_Diagram} are intact, then  $U(K) \leq (0,1)$.   In particular, if $K$ is nontrivial, then $U(K)=(0,1)$.  

\vspace{0.3cm}

Let $K$ be a virtual knot presented by a diagram $D$ which is obtained from a standard 
diagram $D_0$ of a twisted knot by virtualizing some crossings such that exactly one of the crossings 
labeled 1 and 2  ( in Fig.~\ref{Twisted_Knot_Diagram}) is virtualized. 
The Gauss diagram $G$ corresponding to $D$ is as shown in Fig.~\ref{c_chord}, where the directions of horizontal chords are shown as an example.  (The directions of horizontal chords depend on the parity of the number of crossings of $D_0$ and on the positions where we apply virtualization.) 
For convenience, let $c$ be the chord which intersects all other chords in $G$ as in Fig.~\ref{c_chord}.  Let $l$ (or $r$) be the number of chords of $G$ intersecting with $c$ 
from left to right (or right to left) when we see the Gauss diagram as in Fig.~\ref{c_chord}, where we forget the direction of $c$. Recall that all horizontal chords have positive writhe.

\begin{exmp}
\label{unknotting_number_twisted}
Let $K$ be a virtual knot presented by a diagram $D$ which is obtained from a standard 
twisted  knot diagram by virtualizing some crossings such that exactly one of the crossings 
labeled 1 and 2  ( in Fig.~\ref{Twisted_Knot_Diagram}) is virtualized.  Then 
\[U(K)=
 \begin{cases}
  (0,  | J(K) / 2 |),      & \mbox{ if $ | l - r | \leq 1$ } \\
  (1,0),              & \mbox{otherwise}   
   \end{cases}
\]
where $l$ and $r$ are the numbers described above.
\end{exmp}

\noindent \textbf{Proof:} 
Let $G$ be the Gauss diagram of $D$.

\begin{figure}[!ht]
\begin{center}
\includegraphics[scale=.5]{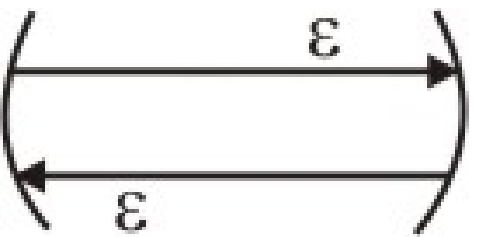}
\caption{ }
\label{fig:t2}
\end{center}
\end{figure}

1. Suppose that $ |l - r | \leq 1$.   We observe that by crossing changes, i.e., changing the directions and the signs, of some chords in $G$, we can change $G$ into a Gauss diagram presenting the trivial knot.  If there exists a pair of chords as in Fig.~\ref{fig:t2}  in $G$, then by a crossing change and an RII move, we can remove the pair. If the resulting Gauss diagram still has such a pair of chords, then we remove the chords by a crossing change and an RII move. Repeat this procedure until we get a Gauss diagram $G_1$ whose horizontal chords are all directed in the same direction. 

\begin{itemize}
 \item
When $ l=r $, the Gauss diagram $G_1$ has only one chord $c$, which presents a trivial knot.  The number of crossing changes is $l$.  On the other hand, all chords of $G$ except $c$ 
are odd chords having positive writhe. The odd writhe number $J(K)$  is
$2l$. By Corollary~\ref{boundsUbyA},  $U(K) = (0, l) = (0,  | J(K) / 2 |)$.  
\item 
When $l = r+1$, the Gauss diagram $G_1$ has two chords, $c$ and a horizontal chord $c_1$, by $r$ times of crossing changes.  Since all horizontal chords are positive, $J(K) = J(G) =  l + r + w(c)$.
\begin{itemize}
\item 
Suppose that $w(c)=-1$.  In this case, the Gauss diagram $G_1$ presents the trivial knot.  Thus 
we have $ U(K) \leq (0,r)$.  On the other hand, $J(K) = 2r$.  
By Corollary~\ref{boundsUbyA}, we have $ U(K) =(0, | J(K)/2 | )$.  
\item 
Suppose that $w(c)=1$.  Apply a crossing change to $G_1$ at $c_1$  and obtain a Gauss diagram $G_2$ with two chords $c$ and $c_2$ where $c_2$ is the horizontal chord obtained from $c_1$ with $w(c_2)=-1$.  Since $G_2$ presents the trivial knot, we have $ U(K) \leq (0,l)$.  On the other hand, $J(K) = 2l$.  
By Corollary~\ref{boundsUbyA}, we have $ U(K) =(0, | J(K)/2 | )$.  
\end{itemize}
\item 
When $r = l+1$, by a similar argument as above, we see that $ U(K) =(0, | J(K)/2 | )$.  
\end{itemize}

2.  Suppose that $ |l - r | > 1$.  
Let $c_{1},c_{2} \dots, c_{n}$ be the horizontal chords of $G$.  For $i=1, \dots, n$, the index ${\rm Ind}(c_i)$ is $1$ or $-1$.    
Note that ${\rm Ind}(c) = r - l $ or $ l - r$ according to the direction of $c$ is upward or downward.  
Let $k = {\rm Ind}(c)$. Then $J_k(G)=w(c)$ and $J_{-k}(G)=0$.  By Proposition~\ref{prop:boundsUbyS}, 
$ (1,0) \leq U(K)$.  

Removing the chord $c$ from $G$, we obtain a Gauss diagram presenting the trivial knot.  
Hence, $U(K) = (1,0)$. \hfill $ \square $

 \vspace{0.3cm}

We used writhe invariants for lower bounds for the unknotting number. A research using the arrow polynomial 
(Dye--Kauffman~\cite{dye2009virtual}) or equivalently Miyazawa polynomial (Miyazawa~\cite{miyazawa2008, miyazawa2009}) 
by the first and the fourth authors will be discussed elsewhere. 

\section{From the virtual knot module}
\label{sect:main_module}

In this section we discuss the unknotting index from a point view of the virtual knot module. 

The group of a virtual knot is defined by a Wirtinger presentation
obtained from a diagram as usual as in classical knot theory 
 \cite{kauffman1998virtual}. 
A geometric interpretation of the virtual knot group is given in \cite{KamKam} and 
a characterization of the group is in \cite{Satoh2000, SW}.   

For the group $G$ of a virtual knot $K$,  let $G'$ and $G''$ be the commutator subgroup 
and the second commutator subgroup of $G$.  Then  the quotient group $M=G'/G''$ forms a finitely generated 
$\Lambda$-module, called the {\it virtual knot module}, where  $\Lambda= {\mathbb Z}[t,t^{-1}]$ denotes 
the Laurent polynomial  ring.  A characterization  of a virtual knot module is given  by 
\cite[Theorem~4.3]{Kaw}.  Let $e(M)$ be the minimal number of $\Lambda$-generators of $M$.   
Then the  following. 

\begin{thm}\label{thm:module}

If $U(K)=(m,n)$, then $e(M)\leq m + n$. 

\end{thm}

Theorem~\ref{thm:module} is obtained from the following lemma. 

\begin{lem}\label{lem:moduleB}
Let $K$ and $K_1$ be virtual knots, and let $M$ and $M_1$ be their virtual knot modules. 
Suppose that a diagram of $K_1$ is obtained from a diagram of $K$ by $(1)$ 
virtualizing a crossing or 
by $(2)$ a crossing change.  Then  $|e(M)-e(M_1)|\leq 1$.
\end{lem} 

\noindent \textbf{Proof:} 
Suppose that a diagram $D_1$ of $K_1$ is obtained from a diagram $D$ of $K$ by 
virtualizing a crossing or a crossing change at a crossing $c$.   

Let $P=(x_1,x_2,\dots, x_u | r_1, r_2,  \dots, r_v)$ be a Wirtinger presentation of  the group $G$ of $K$ 
obtained from $D$ 
with edge generators $x_1,x_2,\dots, x_u$ such that  the last word $r_v$ is $x_u^{-1} x_i^{-1}x_{u-1}x_i$ for an $i<u-1$ which is a relation around the 
crossing $c$.  
  
Let $P_0=(x_1,  x_2,\dots, x_{u-1} |\, \tilde r_1,  \tilde r_2,  \dots, \tilde r_v)$ be the group presentation of a group $G_0$ 
obtained from $P$ by writing the letter $x_u$ to $x_{u-1}$ and rewriting the words $r_1, r_2,  \dots, r_v$ as the words 
$\tilde r_1,  \tilde r_2,  \dots, \tilde r_v$ in  the letters $x_1,  x_2,\dots, x_{u-1}$. Then  the word $\tilde r_v$ is given by 
$x_{u-1}^{-1} x_i^{-1}x_{u-1}x_i$.  

(1) We consider the case of virtualizing the crossing $c$.  
The group $G_1$  of $K_1$ has  the  Wirtinger presentation 
$P_1=(x_1,  x_2,\dots, x_{u-1} |\, \tilde r_1,  \tilde r_2,  \dots, \tilde r_{v-1})$. 

We have  the  following $\Lambda$-semi-exact sequence 
\[ \Lambda[r_1^*, r_2^*,  \dots, r_v^*] \stackrel{d_2}{\to} \Lambda[x_1^*,x_2^*,  \dots, x_u^*] \stackrel{d_1}{\to} \Lambda
\stackrel{\varepsilon}{\to} {\mathbb Z} \to 0\]
of the presentation $P$ by  using the fundamental formula of the Fox differential calculus  in \cite{CroFox}, 
where $\Lambda[r_1^*, r_2^*,  \dots, r_v^*]$ and $\Lambda[x_1^*,x_2^*,  \dots, x_u^*]$ are  free $\Lambda$-modules 
with bases $r_i^*$ ($i=1,2,\dots, v$) and $x_j^*$ ($j=1,2,\dots, u$),  respectively, and the $\Lambda$-homomorphisms 
$\varepsilon$, $d_1$  and $d_2$ are given as follows:
\[\varepsilon(t)=1,\,\, d_1(x_j^*)=t-1\,(j=1,2,\dots, u),\,\, d_2(r_i^*)=\sum_{j=1}^u \frac{\partial r_i}{\partial x_j} x^*_j\, (i=1,2,\dots, v) \]
for the Fox differential calculus $\frac{\partial r_i}{\partial x_j}$ regarded as an element of $\Lambda$ by letting $x_j$ to $t$.  
(Here a $\Lambda$-semi-exact sequence means that, in the above sequence, it is a chain complex of $\Lambda$-modules with $\mbox{Im}(d_1) = \mbox{Ker}(\varepsilon)$ and $\mbox{Im}(\varepsilon) = {\mathbb Z}$.) 
The $\Lambda$-module $M$ of $G$ is  identified with the quotient $\Lambda$-module $\mbox{Ker}(d_1)/\mbox{Im}(d_2)$.

Similarly, we have  the  $\Lambda$-semi-exact sequences 
\begin{eqnarray*}
\Lambda[\tilde r_1^*,\tilde  r_2^*,  \dots, \tilde r_{v-1}^*, \tilde r_v^*] &\stackrel{\tilde d_2}{\to} &
\Lambda[x_1^*, x_2^*,  \dots, x_{u-1}^*] \stackrel{\tilde d_1}{\to} \Lambda
 \stackrel{\varepsilon}{\to} {\mathbb Z} \to 0, \quad \mbox{and} \\
\Lambda[\tilde r_1^*,\tilde  r_2^*,  \dots, \tilde r_{v-1}^*] &\stackrel{\tilde d'_2}{\to}& 
\Lambda[x_1^*, x_2^*,  \dots, x_{u-1}^*] \stackrel{\tilde d_1}{\to} \Lambda
\stackrel{\varepsilon}{\to} {\mathbb Z} \to 0
\end{eqnarray*}
of the presentations $P_0$ and $P_1$, respectively,  where the $\Lambda$-homomorphism $\tilde d'_2$ is the restriction of the $\Lambda$-homomorphism $\tilde d_2$. 
The $\Lambda$-modules $M_0$ and $M_1$ of $G_0$ and $G_1$ are identified with the quotient $\Lambda$-modules 
$\mbox{Ker}(\tilde d_1)/\mbox{Im}(\tilde d_2)$ and $\mbox{Ker}(\tilde d_1)/\mbox{Im}(\tilde d'_2)$, respectively. 

The epimorphism $G\to G_0$ induces a commutative ladder diagram from the $\Lambda$-semi-exact sequence of $P$ to the $\Lambda$-semi-exact sequence of $P_0$ 
by sending $ r_i^*$ to $\tilde r_i^*$  for all $i$,  $x_j^*$ to $x_j^*$ for all $j\leq u-1$  and $x_u^*$   to  $x_{u-1}^*$. 
Then the short $\Lambda$-exact sequence 
\[0\to \Lambda [x_u^*-x_{u-1}^*] \to \Lambda [x_1^*,x_2^*,  \dots, x_u^*] \to
\Lambda[x_1^*,x_2^*,  \dots, x_{u-1}^*]\to 0\]
induces a $\Lambda$-exact sequence 
$\Lambda\to M\to M_0\to 0$, giving 
\[e(M_0)\le e(M)\leq e(M_0)+1. \]
On the other hand, the epimorphism $G_1\to G_0$ induces a commutative ladder diagram from the $\Lambda$-semi-exact sequence of $P_1$ to the $\Lambda$-semi-exact sequence of $P_0$ 
by sending $\tilde r_i^*$ to $\tilde r_i^*$  for all $i\le v-1$ and  $x_j^*$ to $x_j^*$ for all $j\leq u-1$. 
Then from the short exact sequence 
\[0\to\mbox{Im}(\tilde d_2)/\mbox{Im}(\tilde d'_2)\to \mbox{Ker}(\tilde d_1)/\mbox{Im}(\tilde d'_2)\to \mbox{Ker}(\tilde d_1)/\mbox{Im}(\tilde d_2)\to 0\] 
and an epimorphism $\Lambda[\tilde r_v^*]\to\mbox{Im}(\tilde d_2)/\mbox{Im}(\tilde d'_2)$, a $\Lambda$-exact sequence 
$\Lambda\to M_1\to M_0\to 0$ is obtained, giving \[e(M_0)\le e(M_1)\leq e(M_0)+1.\]  
Thus,  the inequality $|e(M)-e(M_1)|\leq 1$ is obtained.   

(2) We consider the case of a crossing change at $c$.  
As seen in (1), we have 
\[e(M_0)\le e(M)\leq e(M_0)+1. \] 
Note that the module $M_0$ is the $\Lambda$-module of the group $G_0$, which is obtained from the virtual knot diagram $D$ with $c$ virtualized by adding a relation that two edge generators around the virtualized $c$ commute. 
When we apply the same argument with $D_1$ instead of $D$, 
we have the same module $M_0$ and 
\[e(M_0)\le e(M_1)\leq e(M_0)+1. \]
Thus the inequality $|e(M)-e(M_1)|\leq 1$ is obtained.  
\hfill $ \square $ 

\begin{rmk}{\rm 
We have an alternative and somewhat geometric proof of 
the case of a crossing change in Lemma~\ref{lem:moduleB} as follow:  
Suppose that a diagram $D_1$ of $K_1$ is obtained from a diagram $D$ of $K$ by a crossing change. 
The virtual knot group $G$ is considered as the fundamental group $\pi_1(E,*)$ of the complement $E=X^*\setminus K$ of a knot $K$ 
 in  a singular 3-manifold $X^*$  which is obtained from the product $X=F\times [0,1]$ with $F$ a closed oriented surface  by shrinking $X_0=F\times 0$ 
to a point $*$, where the knot $K$ is in the interior of $X$ (see \cite{KamKam}).  
The virtual knot group $G_1$  is  the fundamental group $\pi_1(E_1,*)$ of a singular 3-manifold $E_1$ obtained from $E$ 
by surgery along a pair of 2-handles.  
Since the virtual knot modules $M$ and $M_1$ are $\Lambda$-isomorphic to the first homology $\Lambda$-modules 
$H_1(\tilde E; {\mathbb Z})$  and $ H_1(\tilde E_1; {\mathbb Z})$ for the infinite cyclic covering 
spaces $\tilde E$ and $\tilde E_1$ of $E$ and $E_1$, respectively, the argument of \cite[Theorem~2.3]{Ka} implies that  
$|e(M)-e(M_1)|\leq 1$. 
}\end{rmk}

\begin{prop}\label{thm:U0n}
For any non-negative integer $n$, there exists a virtual knot $K$ with $U(K)=(0,n)$. 
\end{prop}

\noindent \textbf{Proof:} 
If $n=1$, then let $K$ be the unknot. Assume $n \geq 1$.  
Let $K$ be the $n$-fold connected sum of a trefoil knot without virtual crossings.  Then $e(M)=n$. 
Since we know $U(K) \leq (0,n)$, Theorem~\ref{thm:module} implies $U(K)=(0,n)$. \hfill $ \square $ 

\begin{thm}\label{thm:U1n}
For any non-negative integer $n$, there exists a virtual knot $K$ with $U(K)=(1,n)$. 
\end{thm}

\noindent \textbf{Proof:} 
Let $D_0$ be a virtual knot diagram obtained from the diagram in Fig.~\ref{kishhh} by applying crossing change at the two crossings on the left side, and let $K_0$ be the virtual knot presented by $D_0$. 
The group of $K_0$ is isomorphic to the group of a trefoil knot.  Let $K$ be the connected sum of $K_0$ and 
the $n$-fold connected sum of a trefoil knot such that a diagram $D$ of $K$ is a connected sum of $D_0$ and 
a diagram of the $n$-fold connected sum of a trefoil knot without virtual crossings.  Then the group of $K$ is isomorphic to the group of the $(n+1)$-fold connected sum of a trefoil knot and $e(M)=n+1$.  
Since $U(K) \leq U(D) \leq (1, n)$, by Theorem~\ref{thm:module}, we have $U(K)=(1,n)$ or $(0, n+1)$.  
On the other hand, when we consider the flat projection $\overline{K}$ of $K$, $\overline{K}$ is the same with 
the flat projection $\overline{K_0}$ of $K_0$, which is non-trivial. By Lemma~\ref{lem:flat}, we have $U(K)=(1,n)$. 
 \hfill $ \square $ 

 \vspace{0.3cm}
 
We conclude with a problem for future work.  

Find a virtual knot $K$ with $U(K)=(m,n)$ for a given $(m,n)$.  A  
connected sum of $m$ copies of Kishino's knot (or $K_0$ in the proof of Theorem~\ref{thm:U1n}) 
and $n$ copies of a trefoil knot seems a candidate.

\hfill
\begin{minipage}{0.5\textwidth}
Kirandeep Kaur\\
Department of Mathematics\\
Indian Institute of Technology Ropar,\\
Punjab, India 140001\\
{\it e-mail: kirandeep.kaur@iitrpr.ac.in\/}\\

Seiichi Kamada\\
Department of Mathematics\\
Osaka City University,\\
 Osaka 558-8585, Japan\\
{\it e-mail: skamada@sci.osaka-cu.ac.jp\/}\\  

Akio Kawauchi\\
OCAMI, Osaka City University,\\
 Osaka 558-8585, Japan\\
{\it e-mail: kawauchi@sci.osaka-cu.ac.jp\/}\\

Madeti Prabhakar\\
Department of Mathematics\\
Indian Institute of Technology Ropar,\\
Punjab, India 140001\\
{\it e-mail: prabhakar@iitrpr.ac.in\/}
\end{minipage}


\begin{thebibliography}{99}

\bibitem{BartholomewFenn} 
A. Bartholomew and R. Fenn: 
{\it Quaternionic invariants of virtual knots and links\/}, 
J. Knot Theory Ramifications 17 (2008),  231--251. 

\bibitem{cheng2013polynomial}  
Z. Cheng and H. Gao:  
{\it A polynomial invariant of virtual links\/},  
J. Knot Theory Ramifications 22 (2013), 1341002. 

\bibitem{cheng2014polynomial}  
Z. Cheng: 
{\it A polynomial invariant of virtual knots\/},  
Proc. Amer. Math. Soc. 142 (2014), 713--725. 

\bibitem{CroFox}  R. H. Crowell and R.  H. Fox: Introduction to Knot Theory, Ginn and Co., 
Boston, Mass., 1963. 

\bibitem{Dye2013}
H. Dye: 
{\it Vassiliev invariants from parity mappings\/}, 
J. Knot Theory Ramifications 22 (2013), 1340008.
 
\bibitem{dye2009virtual}  
H. A. Dye and  L. H. Kauffman: 
{\it Virtual crossing number and the arrow polynomial\/}, 
J. Knot Theory Ramifications 18 (2009),  1335--1357.

\bibitem{FennTuraev}
R. Fenn and V. Turaev: 
{\it Weyl algebras and knots}, 
J. Geom. Phys. 57 (2007), 1313--1324.


\bibitem{ImKimLee2013}
Y. H. Im, S. Kim, and D. S. Lee:
{\it The parity writhe polynomials for virtual knots and flat virtual knots\/}, 
J. Knot Theory Ramifications  22 (2013), 1250133.

\bibitem{Kadokami}
T. Kadokami: 
{\it Detecting non-triviality of virtual links\/}, 
J. Knot Theory Ramifications 12 (2003),  781--803.  

\bibitem{KamKam} N. Kamada and S. Kamada: Abstract link diagrams and virtual knots,  
J.  Knot Theory  Ramifications 9 (2000),  93--106.


\bibitem{kauffman1998virtual} 
L. H. Kauffman:
{\it Virtual knot theory\/}, 
European J. Combin. 20 (1999), 663--690.  

\bibitem{kauffman2004} 
L. H. Kauffman:  
{\it A self-linking invariant of virtual knots\/}, 
Fund. Math. 184 (2004), 135--158. 

\bibitem{kauffman2013} 
L. H. Kauffman: 
{\it An affine index polynomial invariant of virtual knots\/}, 
J. Knot Theory Ramifications 22 (2013), 1340007.

\bibitem{Ka} A. Kawauchi:  Distance between links by zero-linking twists, 
Kobe J. Math. 13 (1996), 183--190.

\bibitem{Kaw} A. Kawauchi: The first Alexander Z(Z)-modules of surface-links 
and of virtual links, 
Geometry and Topology Monographs 14 (2008), 353--371.

 
\bibitem{KishinoSatoh}
T. Kishino and S. Satoh:
{\it A note on classical knot polynomials\/},
J. Knot Theory Ramifications  13 (2004),  845--856.


\bibitem{miyazawa2008} 
Y. Miyazawa: 
{\it A multi-variable polynomial invariant for virtual knots and links\/}, 
J. Knot Theory Ramifications 17 (2008),  1311--1326.

\bibitem{miyazawa2009} 
Y. Miyazawa:
{\it A virtual link polynomial and the crossing number\/},
J. Knot Theory Ramifications 18 (2009),  605--623. 

\bibitem{Satoh2000}
S. Satoh: 
{\it Virtual knot presentation of ribbon torus-knots}, 
J. Knot Theory Ramifications 9 (2000), 531--542.

\bibitem{SatohTaniguchi2014}
S. Satoh and K. Taniguchi:
{\it The writhes of a virtual knot\/}, 
Fund. Math. 225 (2014),  327--342.

\bibitem{SW}
D. S. Silver and S. G. Williams:
{\it Virtual knot groups}, 
Knots in Hellas '98 (Delphi), 440--451, Ser. Knots Everything, 24, World Sci. Publ., River Edge, NJ, 2000.

\end{thebibliography}
\end{document}